\title{
More on the bipartite decomposition of random graphs}
\author{Noga Alon
\thanks{Sackler School of Mathematics and
Blavatnik School of Computer Science,
Tel Aviv University,
Tel Aviv 69978, Israel and School of Mathematics,
Institute for Advanced Study, Princeton,
NJ 08540. Email: nogaa@tau.ac.il.
Research supported in part by a
USA-Israeli
BSF grant, by an ISF grant,
by the Israeli I-Core program
and by the Oswald Veblen Fund.}
\and
Tom Bohman
\thanks{Department of Mathematical Sciences, Carnegie Mellon University, Pittsburgh, PA 15213. Email: tbohman@math.cmu.edu. Research supported in part by NSF grant DMS-1362785}
\and Hao Huang
\thanks{Institute for Mathematics and its Applications, Minneapolis, MN 55455. Email: huanghao@ima.umn.edu.}
}
\documentclass[11pt]{article}
\usepackage{amsfonts,amssymb,amsmath,latexsym}
\def\qed{\ifvmode\mbox{ }\else\unskip\fi\hskip 1em plus 10fill$\Box$}

\newtheorem{theo}{Theorem}[section]

\newtheorem{claim}[theo]{Claim}

\newcommand{\FF}{{\cal F}}

\begin{document}
\date{}

\maketitle

\begin{abstract}
For a graph $G=(V,E)$, let $bc(G)$ denote the minimum number of
pairwise edge disjoint
complete bipartite subgraphs of $G$ so that each edge of $G$
belongs to exactly one of them. It is easy to see that for every
graph $G$,
$bc(G) \leq n -\alpha(G)$, where $\alpha(G)$ is the maximum size
of an independent set of $G$. Erd\H{o}s conjectured in the 80s
that for almost every graph $G$ equality holds,
i.e., that for the random graph $G(n,0.5)$,
$bc(G)=n-\alpha(G)$ with high probability,
that is, with probability that tends to $1$ as $n$
tends to infinity. The first author showed that this is slightly
false, proving that
for most values of $n$ tending to infinity
and for $G=G(n,0.5)$,
$bc(G) \leq n-\alpha(G)-1$ with high probability.
We prove a stronger bound: there exists an absolute constant $c>0$
so that $bc(G) \leq n-(1+c)\alpha(G)$ with high probability.

\noindent
\end{abstract}
\section{Introduction}

For a graph $G=(V,E)$, let $bc(G)$ denote the minimum number of
pairwise edge disjoint
complete bipartite subgraphs of $G$ (bicliques of $G$)
so that each edge of $G$
belongs to exactly one of them. A well known theorem of Graham and
Pollak \cite{GP} asserts that $bc(K_n)=n-1$, see \cite{Tv},
\cite{Pe}, \cite{Vi} for more proofs, and \cite{Al}, \cite{KRW} for
several variants.

Let $\alpha(G)$ denote the maximum size of an independent set of
$G$. It is easy to see that for every graph $G$,
$bc(G) \leq n -\alpha(G)$. Indeed one can partition all edges
of $G$ into $n-\alpha(G)$ stars centered at the vertices of the
complement of a maximum independent set in $G$.
Erd\H{o}s conjectured (see
\cite{KRW}) that for almost every graph $G$ equality holds,
i.e., that for the random graph $G(n,0.5)$,
$bc(G)=n-\alpha(G)$ with high probability
({\em whp}, for short),
that is, with probability that tends to $1$ as $n$
tends to infinity.

Chung and Peng \cite{CP} extended the conjecture for the random
graphs $G(n,p)$ with $p \leq 0.5$, conjecturing that for any
$p \leq 0.5,~~ $
$bc(G) =n -(1+o(1)) \alpha(G)$ whp. They also established lower
bounds supporting this conjecture, and the one of Erd\H{o}s,
by proving that for $G=G(n,p)$ and
for all $0.5 \geq p
\geq \Omega(1)$, $bc(G) \geq n-o((\log n)^{3+\epsilon})$ for
any positive  $\epsilon$.

The first author proved in \cite{Al1} that
Erd\H{o}s' conjecture for $G=G(n,0.5)$
is (slightly) incorrect. It turns out that for most values of $n$,
and for $G=G(n,0.5),~~~$
$bc(G) \leq n - \alpha(G)-1$ whp, while for some exceptional
values of $n$ (that is, those values for which the size of
$\alpha(G)$ is concentrated in two points, and not in one),
$bc(G) \leq n-\alpha(G)-2$ with probability that is bounded  away
from $0$.

He also improved the estimates of \cite{CP} for $G(n,p)$ for any
$c \geq p \geq \frac{2}{n}$,
where $c$ is some small positive absolute constant,
proving that if
$\frac{2}{n} \leq p \leq c$ then for $G=G(n,p)$
$$
bc(G)= n-\Theta(\frac{\log (np)}{p})
$$
whp.

In this note we establish a better upper bound for $bc(G)$ for
$G=G(n,0.5)$, as follows.
\begin{theo}
\label{t11}
There exists an absolute constant $c>0$ so that for $G=G(n,0.5)$,
$$
bc(G) \leq n-(2+2c) \log_2 n \leq n-(1+c)\alpha(G)
$$
with high probability.
\end{theo}

The proof is based on an application of
the second moment method applied to an appropriately defined random
variable. We also describe
another argument,
based on a three-stage exposure of the edges of the
random graph, which provides a simple proof of the fact that
for $G=G(n,0.5)$,
\begin{equation}
\label{e11}
bc(G) \leq n-\alpha(G) - \Omega(\log \log n).
\end{equation}
Although this is weaker than the assertion of Theorem
\ref{t11} we believe this proof is also interesting.

The rest of this note is organized as follows. In Section
2 we describe the short proof of (\ref{e11}). Section 3 includes
the proof of Theorem \ref{t11}. The final Section 4 contains some
concluding remarks, open problems and a brief discussion of
related questions.

Throughout the rest of the note we assume, whenever this is needed, that
$n$ is sufficiently large. To simplify the presentation we omit all
floor and ceiling signs whenever these are not crucial. We make no
attempt to optimize the absolute constants in our estimates.

\section{Three stage exposure and the birthday paradox}

In this section we give a proof of
inequality (\ref{e11}) based
on the following two facts:
\begin{itemize}
\item[(1)] If $ p =c$ where $c$ is a constant
then $\alpha(G(n,p)) = 2 \log_bn - 2 \log_b \log_b n + \Theta(1) $
where $ b = 1/(1-p) $ with high probability, and
\item[(2)] If we choose $a$ items uniformly and independently at random
from a collection of $b$ items (with replacement)
then the probability that the $a$ items are
all distinct is at most $ e^{ -a(a-1)/2b}$.
\end{itemize}
The second fact is known as the birthday paradox.

Let $X,Y$ be an equi-partition of the
vertex set of $G$.  We expose the random edges in
three stages: We first observe
edges inside $X$, then we expose the edges
between $X$ and $Y$, and
finally we reveal the edges within $Y$.
It follows from fact (1) that {\em whp} $X$ contains an
independent set $I$ such that
\[ |I| \ge 2 \log_2n - 2 \log_2 \log_2n - O(1). \]
Let $ \ell= (\log_2n)^{1/3}$.
We partition $Y$ into sets $ Y_1, Y_2, \dots, Y_\ell$
of size $ n/( 2 \log_2 n)^{1/3} $.  Note that for
every vertex $v \in Y$ the neighborhood of $v$ in $I$
is a uniform random subset of $ I$.  Thus it follows from
fact (2) that the probability that every vertex
in $Y_i$ has a different neighborhood in $I$ is at most
\[ \exp \left\{ - \Omega \left( \frac{    n^2
/ (\log_2 n)^{2/3} }{ n^2/ (\log_2 n)^2 } \right) \right\}
= \exp\left\{ - \Omega (\log_2n)^{4/3} \right\} = o(1/n).\]
It follows that with high probability
each set $Y_i$ contains a pair $ a_i,b_i $ of distinct vertices
that have the same neighborhood
in $I$.  Let $I_i = I \cap N(a_i) = I \cap N(b_i)$.
Once this collection of pairs is fixed,
we reveal the edges within $Y$.  With high probability
at least $ \ell/3$ of the pairs $a_i,b_i$
are non-edges, and it follows from fact (1) (taking $p =1/16$)
that among these $ \ell/3$ there is a
collection of  $ \Omega( \log \log n ) $ pairs $ a_i, b_i$ that
spans no edge.  We decompose the edge
set of $G$ into $ n - |I| - \Omega( \log\log n) $ bicliques using
the bicliques $ \{a_i, b_i\} \times I_i $
for the pairs $ a_i,b_i$ in this collection together
with a collection of stars.
\hfill$\Box$\medskip

\section{The proof of the main result}

The proof of Theorem \ref{t11} is based on the second moment method.
The crucial point here is the choice of the random variable to which
it is applied.

For a (large) integer $k$ define a family $\FF_k$ of graphs on
$k$ vertices, as follows.  Each graph in $\FF_k$ is a bipartite
graph with classes of vertices $A$ and $B$, where $|A|=0.1k$ and
$|B|=0.9k$. The set $A$ is the disjoint union of $r=0.01k$
sets  $A_1, A_2, \ldots ,A_r$, where $|A_i|=10$ for each $i$.
For each vertex $b \in B$ there is a binary vector $v_b=(v_b(1),v_b(2),
\ldots ,v_b(r))$ of length $r$. If $v_b(i)=0$ then there are no
edges between $b$ and $A_i$, and if $v_b(i)=1$ then $b$ is
connected to all members of $A_i$. We further assume (although this
is not too crucial, but simplifies matters) that all the vectors
$\{v_b: b \in B\}$ are distinct,
and that the degree of each
$a \in A$ is at least $k/3$ (that if, for each $i$, $v_b(i)=1$ for
at least $k/3$ indices $i$.) In addition we assume that for each
two distinct $i,j$ corresponding to
different sets $A_i, A_j$, the number of
vertices $b \in B$ so that $v_b(i) \neq v_b(j)$ is at least $k/3$.
The family $\FF_k$ contains all the
above graphs.

Note that each graph $F \in \FF_k$ is a bipartite graph on
$k$ vertices satisfying $bc(F) \leq r$. Indeed, the $r$ complete
bipartite graphs with classes of vertices $A_i$ and
$\{b \in B: v_b(i)=1\}$, $(1
\leq i \leq r)$ form a bipartite decomposition of $F$. Let $f_k$
denote the number of graphs on $k$ labelled vertices that are
members of $\FF_k$. We claim that
\begin{equation}
\label{e31}
f_k =(1-o(1)){k \choose {10}} {{k-10} \choose {10}} \ldots
{{k-10r+10} \choose {10}} \frac{1}{r!}  (2^r)^{0.9k}
=k^{(0.09+o(1))k} 2^{0.9rk}.
\end{equation}
Indeed, there are
$$
{k \choose {10}} {{k-10} \choose {10}}\ldots
{{k-10r+10} \choose {10}} \frac{1}{r!}
$$
ways to choose the disjoint sets $A_1,A_2, \ldots, A_r$. After
these are chosen, there are $2^r$ possibilities to choose the edges
from $b$ to the sets $A_i$, for each of the $0.9k$ vertices of
$B$. For a typical choice of these edges, the degree of each
$a \in \cup A_i$ is close to $0.5 \cdot 0.9k $ with high
probability, no two vertices of $B$ have the same sets
of neighbors, and the symmetric difference between the
sets of neighbors of any two vertices of $A$ belonging to distinct sets
$A_i$ is also close to $0.9k/2$. This means that indeed
$1-o(1)$  of the above choices lead to distinct members of $\FF_k$,
establishing (\ref{e31}). For our purpose here it suffices to
note that by the above, since $r=0.01k$,
$$
f_k=2^{0.9rk} k^{\Theta(k)} =2^{(0.9+o(1)) rk}.
$$

Let $V=\{1,2, \ldots ,n\}$ be a fixed set of $n$ labeled vertices,
and let $G=G(n,0.5)=(V,E)$ be the random graph on $V$.
Let $h(k)={n \choose k}f_k2^{-{k \choose 2}}$ be
the expected number of members of $\FF_k$ that appear as induced
subgraphs of $G$ and (with a slight abuse of notation) let
$k$ be the largest integer such that $h(k) \geq 2^k$.
It is not difficult to check that this value of $k$ satisfies
$$
k=2 \log_2 n + 1.8r +O(\log k) =(1+o(1)) 2 \log_2 n +0.018 k
$$
implying that $k=(1+o(1)) \frac{1}{0.982} 2 \log_2 n$, which is
slightly bigger than
$2.036 \log_2 n$. Note that $r=0.01k <0.0204 \log_2 n$.
If $G$ contains an induced copy of a member $F$ of $\FF_k$ then
$$
bc(G) \leq n-k+bc(F) \leq n-2.036 \log_2 n+0.0204 \log_2 n \leq
n-2.015 \log_2 n.
$$
As it is well known that $\alpha(G)=(2+o(1)) \log_2 n$ whp (see
\cite{BE}, \cite{AS}), it suffices to show that $G$ contains
such an induced subgraph whp in order to complete the proof of the
theorem. We proceed to do so using the second moment method.

For each $K \subset V$, $|K|=k$, let $X_K$ be the indicator random
variable whose value is $1$ iff $K$ induces a member of $\FF_k$ in
$G$.
Let $X=\sum_K X_K$, where $K$ ranges over all subsets of size $k$
of $V$,  be the total number of such induced members.
The expectation of this random variable is
$E(X)=h(k)\geq 2^k$. We proceed to estimate
its variance. For $K,K' \subset V$, $|K|=|K'|=k$, let $K \sim K'$
denote that $|K \cap K'| \geq 2$ (and $K \neq K'$). The variance
of $X$ satisfies:
\begin{equation}
\label{e39}
\mbox{Var}(X)=\sum_K \mbox{Var}(X_K)+\sum_{K \sim K'}
\mbox{Cov}(X_K,X_{K'}) \leq E(X) +\sum_{K \sim K'} E(X_K X_{K'}),
\end{equation}
where $K,K'$ range over all ordered pairs of subsets of size $k$ of
$V$ satisfying $2 \leq |K \cap K'| \leq k-1$.

For each $i$, $2 \leq i \leq k-1$, let $h_i$ denote the
contribution of the pairs with intersection $i$ to the above sum,
that is
$$
h_i= \sum_{|K \cap K'|=i} E(X_K X_{K'}).
$$
Our objective is to show that $\sum_{i=2}^{k-1} h_i =o(h(k)^2).$

We consider
two possible ranges for the parameter $i$, as follows.
\vspace{0.2cm}

\noindent
{\bf Case 1:}\, $2 \leq i \leq 0.9k$.
In this case
$$
h_i \leq {n \choose k} f_k {k \choose i} {{n-k} \choose {k-i}} f_k
2^{-2 {k \choose 2}+{i \choose 2}}.
$$
Indeed, for each of the ${n \choose k} f_k$ choices of the set
$K$ and the induced subgraph on it which is a member of $\FF_k$,
there are ${k \choose i} {{n-k} \choose {k-i}}$ ways to choose the
set of vertices $K'$ and then at most $f_k$ ways to  select the
induced subgraph on $K'$ (here there is an inequality, as many of
these choices could lead to inconsistent assumptions about the
induced subgraph on $K \cap K'$, but in this case we have enough
slack and this trivial inequality suffices).
Therefore
\begin{equation}
\label{e33}
\frac{h_i}{h(k)^2} \leq
\frac{{k \choose i} {{n-k} \choose {k-i}}}{{n
\choose k}}  2^{{i \choose 2}} \leq k^i \left(\frac{k}{n}\right)^i 2^{{i
\choose 2}} =\left(\frac{k^2 2^{(i-1)/2}}{n}\right)^i \leq \frac{1}{n^{0.05i}}.
\end{equation}
Here we used the facts that $k \leq 2.04 \log_2 n$ and
$i \leq 0.9k$ to conclude that
$$
\frac{k^2 2^{(i-1)/2}}{n} < \frac{1}{n^{0.05}}.
$$
\vspace{0.2cm}

\noindent
{\bf Case 2:}\, $i=k-j$ where $1 \leq j \leq 0.1k$. This case is
more complicated and requires a careful estimate of the number of
possibilities for the induced subgraph on $K \cup K'$.This is done
in the following claim.
\begin{claim}
\label{c32}
Let $k$ and $r=0.01k$ be as above, and let
$K,K'$ be two sets of labelled vertices, where
$|K|=|K'|=k$ and $|K \cap K'|=i=k-j$ with $1 \leq j \leq 0.1k$.
Then the number of graphs $H$ on $K \cup K'$, such that the induced
subgraph of $H$ on $K$ and the induced subgraph of $H$ on $K'$ are
members of $\FF_k$ is at most
$$
f_k (r+2^r)^j 2^{0.9kj/10}
$$
\end{claim}
\vspace{0.1cm}

\noindent
{\bf Proof of Claim:}\, There are $f_k$ ways to choose the induced
subgraph of $H$ on $K$. Fixing such a choice, we estimate the
number of ways to extend it to the edges inside $K'$ (which are not
inside $K$, as this part is already fixed). Let $A$ and $B$ denote
the vertex classes of the member $F'$
of $\FF_k$ in $K'$, thus
$A \cup B=K'$. Let $A=A_1 \cup A_2 \ldots \cup A_r$ denote the
partition of $A$ into disjoint sets of size $10$ in this member.
Since $|K \cap K'| \geq 0.9k$ and the degree of each $A$-vertex in
$F'$ is at least $k/3$ whereas the degree of each $B$-vertex is at
most $|A|=0.1k$ it follows that any vertex $a \in A$ must have
at least $k/3-0.1k >0.1k$ neighbors in $K \cap K'$, and thus
knowing the edges inside $K \cap K'$ reveals the fact that this is
an $A$-vertex.  We thus know, for each vertex in $K \cap K'$,
if it is an $A$-vertex or a $B$-vertex. Moreover, since the
sets of $B$-neighbors of any two $A$ vertices from distinct subsets
$A_i$ differ on at least $k/3$ vertices $b \in B$, the edges inside
$K \cap K'$ reveal, for each $i$ so that $A_i$ intersects $K \cap
K'$, all the vertices of $A_i \cap (K \cap K')$.
There are now at most $(r+2^r)^j$ ways to choose, for each vertex
in $K'-K$, if it lies in one of the sets $A_i$ (which is either
represented in $K \cap K'$ or not), and if so, decide to which of
the $r$ sets it belongs, and in addition, if it is a $B$-vertex,
to decide to which sets $A_i$ it is connected.
Here we are over-counting, as we ignore the
fact that any set $A_i$ has to be of size exactly $10$,  but this
estimate suffices.  Note that after the above choices, the identity
of all vertices in each set $A_i$ is known. As each set $A_i$ is of
cardinality $10$, there are at most $j/10$ sets $A_i$ which are
completely contained in $K'-K$. For each such set, there are
at most $2^{0.9k}$ possibilities  to choose the edges between the
vertices of this set and the remaining vertices of $K'$. Once these
choices are made, all edges  inside $K'$ are determined.
This completes the proof of the claim
\hfill$\Box$\medskip

Returning to the proof of the theorem, we proceed with the estimate
of $h_i/h(k)^2$ in Case 2. By the claim, for $i=k-j, j \leq 0.1k$
we have (since $h(k) \geq 2^k>1$):
$$
\frac{h_i}{h(k)^2} \leq \frac{h_i}{h(k)}
\leq {k \choose j} {{n-k} \choose j} (r+2^r)^j 2^{0.9kj/10}
2^{-(k-j)j}
$$
$$
\leq [kn 2^{2r} 2^{0.9k/10}2^{-(k-j)}]^j \leq n^{-0.5j},
$$
with room to spare.

Combining the last inequality with (\ref{e39}) and (\ref{e33}),
and using the fact that $E(X) =h(k) \geq 2^k$,  we conclude that
$Var(X) =o(E(X)^2)$ and hence, by Chebyshev's Inequality,
$X>0$  whp. This implies that $bc(G) \leq n-2.015 \log_2 n$ whp,
completing the proof of Theorem \ref{t11}.
\hfill$\Box$\medskip


\section{Concluding remarks and open problems}
\begin{itemize}
\item The estimate in Theorem~\ref{t11} is the best we can hope
to get with this method, up to the constant $c$. This is because
all members of $\FF_k$ are bipartite graphs, and the random
graph $G=G(n,0.5)$ cannot contain any induced bipartite graph on more
than $2 \alpha(G)$ vertices.
\item
We have shown that
for
$G=G(n,0.5)$, $bc(G) \leq n-\alpha(G)-\Omega(\log n)$
whp. It will be  interesting to decide whether or not
$bc(G)=n-O(\alpha(G))$ whp.
\item
For  $p<0.5$ and $G=G(n,p)$
it seems that both proofs we know  do not give
any improvement of the trivial estimate
$bc(G) \leq n-\alpha(G)$. Is it true that for any
fixed positive $p <0.5$, $bc(G)=n-\alpha(G)$ whp ? (for $p>1/2$ it is
easy to get a better upper bound).
\end{itemize}

We conclude this short
paper with a note regarding biclique decompositions of
twin-free graphs.  Vertices $u$ and $v$ in a
graph $G$ are twins if they have exactly the same
neighborhoods, and $G$ is twin-free if $G$ contains no such pair
of vertices.  Note that if a pair of vertices $u,v$ are twins in $G$ then
$ bc(G) = bc(G-u)$.  Thus it is quite natural to consider the maximum number
of vertices in a twin-free graph $G$ with $ bc(G) = r$.
\begin{theo}
Suppose $G$ is a twin-free graph whose edges can be decomposed into $r$ bicliques, then $|V(G)| \le 2^{r+1}-1$ and this bound is tight.
\end{theo}
{\bf Proof:}
We first construct a graph $G$ which attains this upper bound.
Let $V(G)$ be a collection of vectors $v$ in $\{0,1,2\}^r$, such that
$v_i=1$ for at most one index $i$, and $v_j=2$ for all $j>i$,
$v_j \in \{0, 2\}$ for all $j<i$. In other words,
$$V(G)=\cup_{k=0}^r\{0,2\}^k \times \{1\} \times \{2\}^{r-k-1}.$$

The number of vertices in $G$ is equal to $1+2+\cdots+2^r=2^{r+1}-1$.
We define two vertices $u$ and $v$ to be adjacent if there exists $i$
such that $(u_i,v_i)=(1,0)$ or $(0,1)$. To show that $G$ is twin-free,
suppose $u$ and $v$ are two distinct vertices of $G$. If $u_i=v_i=1$
for some $i$, then one can find $j<i$ so that $(u_j,v_j)=(2,0)$ or $(0,2)$,
then the vector $w$ with $w_j=1$ and $w_k=2$ for all $k \ne j$ is only
adjacent to one of $u$ and $v$. If $u_i=1$ for some $i$ and $v_i \ne 1$, then the vector
$w$ with $w_i=0$ and $w_j=2$ for all $j \ne i$ is adjacent to $u$ but not $v$.
Finally if both $u$ and $v$ are in $\{0,2\}^r$, take the coordinate $i$ such
that $(u_i, v_i)=(0,2)$ or $(2,0)$, then again letting $w_i=1$ and $w_j=2$
for all $j \ne i$ shows that they have different neighborhoods.

The definition of $G$ naturally induces an edge decomposition into bicliques:
two vertices $u$ and $v$ are adjacent in the biclique $G_i$ iff $(u_i, v_i)=(0,1)$
or $(1,0)$. To verify that this is indeed a partition, assume that the edge $uv$
belongs to two bicliques $G_i$ and $G_j$. This can only happen when $(u_i,u_j)=(0,0)$,
$(v_i,v_j)=(1,1)$ or $(u_i,u_j)=(0,1)$, $(v_i,v_j)=(1,0)$ (when necessary we swap $u$
and $v$). Note that both cases are impossible since all the vectors in $V(G)$ have
at most one coordinate equal to $1$, and $0$ never appears after $1$.

Next we are going to show that $2^{r+1}-1$ is an upper bound. For a twin-free graph
$G$ with biclique partition $E(G)=\cup_{i=1}^r E(G_i) = \cup_{i=1}^r E(A_i, B_i)$, we
assign a $r$-dimensional vector $v_u$ to every vertex $u$, such that $(v_u)_i=1$ if
$u \in A_i$, $0$ if $u \in B_i$ and $2$ otherwise. Note that two vertices associated
with the same vector have common neighborhoods, so we may assume that all the vectors
$v_u$ are distinct. Let $\mathcal{F}=\{v_u\}_{u \in G}$, and
$\mathcal{F}_I=\{v: v \in \mathcal{F}, \{i:v_i \in \{0,1\}\}=I \}$. We claim that for
all $|I| \ge 1$, $|\mathcal{F}_I| \le 2$. This is obvious for $|I|=1$. The case
$|I| \ge 2$ follows from the observation that among any three distinct vectors in $\{0,1\}^I$,
there always exists a pair differing in at least two coordinates $i$ and $j$, which
contradicts the assumption that $G_i$ and $G_j$ are disjoint. Therefore
$$|\mathcal{F}| \le 1+\sum_{i=1}^r 2\binom{r}{i}=2^{r+1}-1.$$
\hfill$\Box$\medskip

\noindent
{\bf Acknowledgment}
Part of this work was done during the workshop on
Probabilistic and Extremal Combinatorics
which took place in IMA, Minneapolis in September, 2014.
We would
like to thank IMA and the organizers of the conference
for their hospitality.  We also thank Pat Devlin and
Jeff Kahn for a helpful conversation.

\end{document}